\documentclass[12pt,a4paper]{article}
\usepackage[T2A]{fontenc}
\usepackage[cp1251]{inputenc}
\usepackage[russian]{babel}
\usepackage{amsmath,amssymb,amsthm,amsfonts,amscd}

\begin{document}
\sloppy
\date{}

\title{Geometric approach towards stable homotopy groups of spheres.
The Steenrod-Hopf invariant}
\author{Pyotr M. Akhmet'ev \thanks{This work was supported in part by the London Royal Society
(1998-2000), RFBR 08-01-00663,  INTAS 05-1000008-7805.}}

\maketitle
\newtheorem{theorem}{Теорема}[section]
\newtheorem*{main*}{Основная Теорема}
\newtheorem*{theorem*}{Теорема}
\newtheorem{lemma}[theorem]{Лемма}
\newtheorem{proposition}[theorem]{Предложение}
\newtheorem{corollary}[theorem]{Следствие}
\newtheorem{conjecture}[theorem]{Гипотеза}
\newtheorem{problem}[theorem]{Проблема}

\theoremstyle{definition}
\newtheorem{definition}[theorem]{Определение}
\newtheorem{remark}[theorem]{Замечание}
\newtheorem*{remark*}{Замечание}
\newtheorem*{example*}{Пример}
\def\Z{{\Bbb Z}}
\def\R{{\Bbb R}}
\def\RP{{\Bbb R}\!{\rm P}}
\def\N{{\Bbb N}}
\def\C{{\bf C}}
\def\A{{\bf A}}
\def\D{{\bf D}}
\def\O{{\bf O}}
\def\I{{\bf I}}
\def\E{{\Bbb E}}
\def\fr{{\operatorname{fr}}}
\def\st{{\operatorname{st}}}
\def\mod{{\operatorname{mod}\,}}
\def\cyl{{\operatorname{cyl}}}
\def\dist{{\operatorname{dist}}}
\def\sf{{\operatorname{sf}}}
\def\dim{{\operatorname{dim}}}

 \begin{abstract}

In this paper a geometric approach toward stable homotopy groups
of spheres, based on the Pontrjagin-Thom [P] construction is
proposed. From this approach a new proof of Hopf Invariant One
Theorem by J.F.Adams [A1] for all dimensions except $15,31,63,127$
is obtained.

It is proved that for $n>127$ in the stable homotopy group of
spheres $\Pi_n$ there is no elements with Hopf invariant one. The
new proof is based on geometric topology methods. The
Pontrjagin-Thom Theorem (in the form proposed by R.Wells [W])
about the representation of stable homotopy groups of the real
projective infinite-dimensional space (this groups is mapped onto
2-components of stable homotopy groups of spheres by the
Khan-Priddy Theorem [A2]) by cobordism classes of immersions of
codimension 1 of closed manifolds (generally speaking,
non-orientable) is considered. The Hopf Invariant is expressed as
a characteristic number of the dihedral group for the
self-intersection manifold of an immersed codimension 1 manifold
that represents the given element in the stable homotopy group. In
the new proof the Geometric Control Principle (by M.Gromov)[Gr]
for immersions in a given regular homotopy classes based on
Smale-Hirsch Immersion Theorem [H] is required.
 \end{abstract}

Let $f: M^{n-1} \looparrowright \R^n$, $n= 2^l -1$, be a smooth
immersion of codimension 1.

The characteristic number
$$ \left< w^{n-1}_1(M);[M^{n-1}] \right> = h(f) $$
is called  the stable Hopf invariant or the Steenrod-Hopf
invariant. This characteristic number depends only on the immersed
manifold $M^{n-1}$ itself. The relationship with the definition of
the Steenrod-Hopf invariant in algebraic topology is considered in
[E],[K],[L].

\subsubsection*{Theorem (by J.F.Adams), [A]}
\begin{center}
For $l \ge 4$, $h(f)=0$.
\end{center}

\subsubsection*{Skew-framed immersions}
Let $f: M^{n-k} \looparrowright \R^n$ be an immersion of
codimension $k$. Let $\kappa: E(\kappa) \to M^{n-k}$ be a line
bundle over $M^{n-k}$ and let $\Xi: k \kappa \to \nu(f)$ be an
isomorphism of the normal bundle of the immersion $f$ with the
Whitney sum of $k$ copies of the line bundle $\kappa$.

We shall call the triple $(f, \kappa, \Xi)$ a skew-framed
immersion with characteristic class $\kappa \in
H^1(M^{n-k};\Z/2)$. (If $n$ is odd then $w_1(\kappa)$ is the
orientation class of $M^{n-k}$, see [A-E] for more details).

\subsubsection*{The Steenrod-Hopf invariants for skew-framed immersions}

The characteristic class
$\left<w_1(\kappa)^{n-k};[M^{n-k}]\right>=h(f,\kappa,\Xi)$ is
called the Steenrod-Hopf invariant of the skew-framed immersion
$(f,\kappa,\Xi)$.
\subsubsection*{The Main Theorem 1}

Let ($f: M^{n-k} \looparrowright \R^n, \kappa, \Xi)$ be a
skew-framed immersion, $n=2^l-1$, ${\rm
dim}(M)=n-k=\frac{n+1}{2}+7$. Then for $n \ge 255$ (i.e. for $l
\ge 8$) $$h(f,\kappa, \Xi)=\left<w_1(\kappa)^{{\rm
dim}(M)};[M]\right>=0 \quad ({\rm mod} \quad 2).$$

\subsubsection*{Corollary}
Adams' Theorem for $n \ge 255$.

\[  \]
Let $\D_4$ be the dihedral group of  order $8$,
$$ \D_4 = \{a,b \vert a^4=b^2=e, [a,b]=a^2 \}.$$
This is the group of symmetries of the two coordinate axes in the
plane. Let $$\I_a = \{e, a, a^2, a^3 \}, \I_b =\{e, b, a^2, a^2b
\}, \I_c = \{e, ab, a^2, a^3b \}$$ be the subgroups in $\D_4$ of
index 2. The cyclic $\Z/4$--subgroup $\I_a$ is generated by
rotation of the plane through the angle $\frac{\pi}{2}$. The $\Z/2
\oplus \Z/2$--subgroup $\I_b$ ($\I_c$) is generated by symmetries
(or reflections) with respect to the bisectors of the coordinate
axes. The $\Z/2$--subgroup $\I_b \cap \I_c = \{e, a^2\}$ is
generated by the central symmetry.


Let ($f: M^{n-k} \looparrowright \R^n, \kappa, \Xi)$ be a
skew-framed generic immersion. Its double points manifold
$N^{n-2k}$ is immersed into $\R^n$, $g: N^{n-2k} \looparrowright
\R^n$, and the normal bundle $\nu(g)$ admits a canonical
decomposition $\Psi: \nu(g)=k \eta^{\ast}$, where $\eta^{\ast}$ is
a two-dimensional bundle over $N^{n-2k}$ with $\D_4$-structure.
The bundle $\eta^{\ast}$ is a pull-back of the universal bundle
$E(\D_4) \to K(\D_4,1)$, via the classifying map $\eta: N^{n-2k}
\to K(\D_4,1)$.

Let us consider the canonical 2-fold covering $\bar N^{n-2k} \to
N^{n-2k}$ over the double point manifold of the immersion $g$. This covering
corresponds to the subgroup $\I_c \subset \D_4$,
$\I_c=\{e,a^2,ab,a^3b\}$.

Let $\bar \kappa \in H^1(\bar N^{n-2k};\Z/2)$ be the cohomology
class corresponding to the epimorphism $\I_c \to \I_d$ with the
image $\I_d=\{e, a^2 \simeq ab \}=\Z/2$ (the kernel is generated
by the element $a^3b$). By the definition $\bar \kappa =
i^{\ast}(\kappa)$, $i: \bar N^{n-2k} \looparrowright M^{n-k}$ is
the canonical immersion of the double point covering. Let us
define the following characteristic number $$\bar h(g,\eta,\Psi)=
\left<\bar \kappa^{n-2k};[\bar N^{n-2k}]\right>.$$

\subsubsection*{Lemma 2}
 $$ h(f, \kappa, \Xi)= \bar h(g, \eta, \Psi). $$

\subsubsection*{Proof of Lemma 2}
Immediate from  Herbert's Theorem. (Concerning Herbert's Theorem,
see e.g. [E-G].)

\subsubsection*{Definition (Cyclic structure for skew-framed immersions)}
Let $(f,\kappa, \Xi)$ be a skew-framed immersion, $N^{n-2k}$ be
the (odd-dimensional) double self-intersection point manifold of
$f$.  A mapping $$\mu: N^{n-2k} \to K(\I_a,1)$$ ($\I_a =
\{e,a,a^2,a^3\}$) is called a {\it cyclic structure} for $f$ if
$$\left<\mu^{\ast}(t);[N^{n-2k}]\right>\quad =\quad h(f),$$ where $t \in
H^{n-2k}(K(\I_a,1);\Z/2)$ is the generator.


The following lemma is proved by an explicit calculation.

\subsubsection*{The Main Lemma 3 (jointly with P.J.Eccles (1998))}
Let $n-k=\frac{n+1}{2}+7$ ($n-2k=15$), $n \ge 31$, $\mu: N^{n-2k}
\to K(\D_4,1)$ be a cyclic structure for $f$. Then
$h(f,\kappa,\Xi)=0 \quad ({\rm mod}  2)$.

\subsubsection*{Lemma 4. Cyclic Structure for skew-framed immersions}
For $n \ge 255$, an arbitrary skew-framed immersion $f: M^{n-k}
\looparrowright \R^n$,
 $n-k=\frac{n+1}{2}+7$ is regularly homotopic to an immersion with a cyclic structure.

\subsubsection*{Proof}
The proof is a corollary of Lemma 5 and Proposition 6.


\subsubsection*{Definition (Cyclic Structure for generic mappings of the standard projective space)}
Let $g: \RP^{n-k} \to \R^n$ be a generic mapping, $n \ge 3k$,
$n=2^l-1$ with double point manifold $N^{n-2k}$ and critical
points $(n-2k-1)$-dimensional submanifold $(\partial N)^{n-2k-1} \subset \RP^{n-k}$.

Let $\eta: (N^{n-2k},\partial N) \to (K(\D_4,1),K(\I_b,1))$ be the
structured mapping corresponding to $g$. This structured mapping
is defined analogously with the case of skew-framed immersions.
Obviously, the restriction of the mapping $\eta$ to the boundary
of the double points manifold, i.e. to the critical points
submanifold, has the target $K(\I_b,1) \subset K(\D_4,1)$. The
standard inclusion $\I_d \subset \I_a$ as the subgroup of the
index 2 is well-defined.
 We shall call a mapping
$$ \mu: (N^{n-2k},\partial N^{n-2k-1}) \to
(K(\I_a,1),K(\I_d,1))$$ a cyclic structure for $g$, if the
following conditions hold:

$\qquad (\ast) \qquad $ the homological condition:
$$ \qquad  \left<\mu^{\ast}(t);[N^{n-2k},\partial N^{n-2k-1}]\right>=1 \quad (mod \quad 2),$$
where $t \in H^{n-2k}(K(\I_a,1),K(\I_d,1);\Z/2)$ is the generator
in the cokernel of the homomorphism
$H^{n-2k}(K(\I_a,1),K(\I_d);\Z/2)$ (remark: this relative
characteristic number is well-defined, by explicit calculations
$\mu_{\ast}([\partial N^{n-2k-1}]) \in H_{n-2k-1}(K(\I_d,1;\Z/2))$
is trivial);
\[  \]
 $\qquad (\ast
\ast) \qquad $ the boundary condition :
$$p_c \circ \eta \vert_{\partial N} : \partial N^{n-2k-1} \to
K(\I_b,1) \to K(\I_d,1),$$ $\I_b = \{e, b, a^2, ba^2 \}$, where
$p_c: K(\I_b,1) \to K(\I_d,1)$ is the standard projection with the
image $\I_d= \{e, a^2 \simeq ba^2\}$, coincides with $\mu
\vert_{\partial N^{n-2k-1}}$.

\subsubsection*{Lemma 5.  Geometric Control}
Let $$(f,\kappa,\Xi), \quad f: M^{n-k} \looparrowright \R^n$$ be
a skew-framed immersion. Let us assume (see Proposition 6 below) that there exists a generic mapping
$$g: \RP^{n-k} \to \R^n, \quad n \ge
3k,$$  with a cyclic structure $$ \mu:
(N^{n-2k},\partial N^{n-2k-1}) \to (K(\I_a,1),K(\I_d,1)).$$ Then
there exists a skew-framed immersion $(f',\kappa, \Xi')$ in the
regular homotopy class of $f$ with a cyclic structure.
\[  \]
\subsubsection*{The idea of the proof of Lemma 5.}
Take the mapping $$g \circ \kappa: M^{n-k} \to \RP^{n-k} \to
\R^n.$$ By [Gr], 1.2.2, in the regular homotopy class of $f$ there
exists a generic immersion   $$f': M^{n-k} \looparrowright \R^n$$
 defined as a $C^0$-close generic regular perturbation (
arbitrary small) of a (singular) mapping $$g \circ \kappa.$$

\subsubsection*{A  construction by S.A.Melikhov (2004)} Let us
denote by $J$ the join of $(2^{l-4}+1)=r$ copies of the standard
$\Z/4$- lens space $S^7/i$, ${\rm dim}(J)=2^{l-1}+7=n-k$. There is a $PL$--embedding
$i_J: J \subset \R^n$ for $l \ge 8$.

Let  $p': S^{n-k} \to J$ be the join of $r$ copies of the standard
cover $S^7 \to S^7/i$, $\hat p: S^{n-k}/i \to J$ be the quotient
mapping  of $p'$, $p: \RP^{n-k} \to J$ be the composition of the
standard projection $\pi: \RP^{n-k} \to S^{n-k}/i$ with $\hat p$.
The composition $i_J \circ \hat p: S^{n-k}/i \to J \to \R^n$ is
well-defined. Let $\hat g: S^{n-k}/i \to \R^n$ be an
$\varepsilon$-small generic alternation of the mapping $i_J \circ
\hat p$, $d: \RP^{n-k} \to \R^n$  be defined by a
$\varepsilon_1$-small generic alternation, $\varepsilon_1 <<
\varepsilon$, of the composition $\hat g \circ \pi$.

\subsubsection*{Proposition 6 \footnote{The author was developed
this proof following conversations with Prof. O.Saeki and Dr.
R.R.Sadykov (2006))}} The Melikhov map $d: \RP^{n-k} \to \R^n$ is
equipped with a cyclic structure.

\[  \]
The rest of the paper concerns the proof of this result. This
Proposition with Lemma 5 implies Lemma 4 and Lemma 4 with the Main
Lemma 3 implies Lemma 2 and the main Theorem 1.

\subsubsection*{The beginning of the proof of the Proposition 6}

Let $\Gamma_0$ be the {\it delated product} of the standard
projective space $\RP^{n-k}$ $$\Gamma_0= (\RP^{n-k} \times
\RP^{n-k} \setminus \Delta_{diag})/T',$$ where the quotient is
determined with respect to the free involution
$$T': \RP^{n-k} \times
\RP^{n-k} \setminus \Delta_{diag} \to \RP^{n-k} \times \RP^{n-k}
\setminus \Delta_{diag}$$ $$T'(x,y)=(y,x).$$
 The classifying map
$$\eta_{\Gamma_0} : \Gamma_0 \to K(\D_4,1)$$ is well-defined.
(Note that $\pi_1(\Gamma_0)=\D_4$ and the involution $T'$
corresponds to the subgroup $\I_c \subset \D_4$.)
Let $\Delta_{antidiag} \subset \Gamma_0$ be a subspace, called the
antidiagonal, defined as $$\Delta_{antidiag} = \{(x,y) \in
(\RP^{n-k} \times \RP^{n-k}/T') \vert T(x)=y \},$$ where
$$T:\RP^{n-k} \to \RP^{n-k}$$ is the standard involution on the
covering $$\RP^{n-k} \to S^{n-k}/i.$$
Let us denote by $\Gamma$ the space $$\Gamma_0 \setminus
(U(\Delta_{antidiag}) \cup U(\Delta_{diag})),$$ where
$U(\Delta_{antidiag})$ is a small regular neighborhoods of
$\Delta_{antidiag}$, $U(\Delta_{diag}$ is a small regular
neighborhood of the end of $\Gamma_0$ near the deleted diagonal
$\Delta_{diag}$. (The radius of the regular neighborhoods depends
on a constant $\varepsilon$ of an approximation in the Melikhov
construction.)

The space $\Gamma$ is a manifold with boundary. The involution
$$T: \RP^{n-k} \to \RP^{n-k}$$ induces the free involution
$$T_{\Gamma}: \Gamma \to \Gamma.$$
A polyhedron $$\Sigma_0= \{[(x, y)] \in \Gamma_0, p(x) = p(y)\},
\quad \Sigma_0 \subset \Gamma_0$$ of double points of the mapping
$$p: \RP^{n-k} \to J$$ is called {\it the singular set} or {\it the singular
polyhedron}.

The mapping
$$\eta_{\Sigma_0}: \Sigma_0 \to K(\D_4,1)$$ is well-defined as the
restriction of the mapping $$\eta_{\Gamma_0} \vert_{\Sigma_0}.$$
The subpolyhedron
$$\Sigma_0 \subset \Gamma_0$$ decomposes as $$\Sigma_0 =
\Sigma_{antidiag} \cup K, \quad K \subset \Gamma, $$ where
$$\Sigma_{antidiag}= \Sigma_0 \cap U(\Delta_{antidiag}).$$ The
restriction $\eta_{\Gamma_0} \vert_{K}$ will be denoted by $\eta_K : K \to K(\D_4,1)$.

\subsubsection*{Boundary conditions of $\eta_K$}

The diagonal and antidiagonal boundary components of $K$ will be
denoted by $$Q_{diag} = K \cap \partial U(\Delta_{diag}), \quad
Q_{antidiag} = K \cap \partial U(\Delta_{antidiag}).$$
\[  \]
The restriction $\eta_K \vert_{Q_{antydiag}} : Q_{antidiag} \to
K(\D_4,1)$ is decomposes as $$i_a \circ \eta_{antidiag}:
Q_{antidiag} \to K(\I_a,1) \subset K(\D_4,1).$$
\[  \]
The restriction $\eta_K \vert_{Q_{diag}} : Q_{diag} \to K(\D_4,1)$ is decomposed as
$$i_b \circ \eta_{diag}: Q_{antidiag} \to K(\I_b,1) \subset
K(\D_4,1).$$

\subsubsection*{The resolution space $RK$.}

Let us construct the space $RK$ called {\it  the resolution space}
for $K$. This space is included into the diagram

$$ K(\I_a,1)   \stackrel {\phi}{\longleftarrow}  RK   \stackrel {pr}{\longrightarrow}
K.
$$
\[  \]
Let us denote $pr^{-1}(Q_{diag})$ by $RQ_{diag}$,
$pr^{-1}(Q_{diag})$ by $RQ_{antidiag}$.
The boundary conditions on $Q_{antidiag}$ are:
$$
\begin{array}{ccc}
RQ_{antidiag} & \qquad \stackrel{pr}{\longrightarrow} \qquad & Q_{antidiag} \\
\phi \searrow & & \swarrow \eta_{antidiag} \\
& K(\I_a,1). &
\end{array}
$$
\[  \]
The boundary conditions on $Q_{diag}$ are:
\[  \]
$$
\begin{array}{ccc}
RQ_{diag} & \qquad \stackrel{pr}{\longrightarrow} \qquad & Q_{diag} \\
\phi \downarrow && \downarrow \eta_{diag}\\
 K(\I_d,1) & \stackrel{p_b}{\longleftarrow} & K(\I_b,1). \\
\end{array}
$$
\[  \]
The diagrams above are included into the following diagram:
$$
\begin{array}{ccc}
 K(\I_a,1) &&\\
\uparrow \phi & &\nwarrow \\
RK  & \longleftarrow & RQ_{diag} \cup RQ_{antidiag} \\
 \downarrow  & & \downarrow  \\
 K
& \supset & Q_{diag} \cup Q_{antidiag} \\
   \downarrow  \eta & & \swarrow    \\
     K(\D_4,1).&& \\
\end{array}
$$



Let us consider Melikhov's mapping $d: \RP^{n-k} \to \R^n$ (this
is
 a small generic alternation of the composition $ i
\circ p \circ  \pi: \RP^{n-k} \to S^{n-k}/i \to J \subset \R^n$).
 Let  $N^{n-2k}$ be the double
point manifold (with boundary) of $d$, the embedding $N^{n-2k}
\subset \Gamma_0$ is well-defined.  The manifold $N^{n-2k}$ is
decomposed into two manifolds (with boundary) along the common
component of the boundary $$N^{n-2k} = N_{antidiag} \cup N_d,$$
$$N_{antidiag}= N^{n-2k} \cap U(\Delta_{antidiag}),$$ $$N_d= N^{n-2k}
\cap \Gamma.$$

\subsubsection*{Lemma 7}

There exists a mapping $ res: N_d \to RK$ called the {\it resolution
mapping} that induces a mapping $\mu: N_d \to K(\I_a,1)$ included
into the following diagram:

$$
\begin{array}{ccccccc}
K(\I_a,1) & =& K(\I_a,1)   & &  \\
\uparrow  \phi  &  & \uparrow \mu & & \nwarrow  \\
RK & \stackrel {res}{\longleftarrow}
&N_d & \supset & W_{diag} \cup W_{antidiag} \\
\downarrow \eta & &  \downarrow \eta & & \swarrow  \\
K(\D_4,1) &= &     K(\D_4,1),&& \\
\end{array}
$$
\[  \]
with boundary conditions on $W_{antidiag}$:
\[  \]
$\mu \vert_{W_{antydiag} \subset N_d} = i_a \circ \eta_{antidiag} : W_{antidiag} \longrightarrow
K(\I_a,1) \stackrel{i_a}{\longrightarrow} K(\D_4,1)$,
\[  \]
and with boundary conditions on $W_{diag}$:
\[  \]
$$\mu = i_a \circ p_b \circ \eta_{diag}: W_{diag} \to K(\I_b,1)
\to K(\I_d,1) \to K(\I_a,1).$$
\[  \]

\subsubsection*{Lemma 8}
The mapping
$$\mu_a= \eta \vert_{N_{antidiag}} \cup \mu: N^{n-2k} =
N_{antidiag} \cup N_d \to K(\I_a,1)$$ determines a cyclic
structure for $d$.

\subsubsection*{Proof of Lemma 8}

We have to prove the equality $\ast$ in the Definition of the
Cyclic Structure for generic mappings.



Let us consider the free involution $T_{\Gamma} : \Gamma \to
\Gamma$ and the quotient $\Gamma/T_{\Gamma}$. The fundamental
group $\pi_1(\Gamma/T_{\Gamma})$, denoted by $\E$, is a quadratic
extension of $\D_4$ by means of an element $c \in \E \setminus
\D_4$, $c^2=a^2$. The element $c$ of  order 4 is commutes with all
elements in the subgroup $\D_4 \subset \E$. The following diagram
is well-defined:



$$
\begin{array}{ccccccc}
N_d & \longrightarrow & RK & \longrightarrow & K & \subset & \Gamma \\
\downarrow & &  \downarrow & & \downarrow & & \downarrow \\
N_d/T & \longrightarrow & RK/T & \longrightarrow  & K/T & \subset & G/T\\
\downarrow & &  \downarrow & & \downarrow & & \downarrow \\
K(\I_a,1) & = & K(\I_a,1) & & K(\E,1) & = & K(\E,1). \\
\end{array}
$$


 The composition
$RK \to RK/T \to K(\I_a,1)$ coincides with $\phi$ (The $c \in
\pi_1(RK/T) \setminus \pi_1(RK)$  commutes with all  elements in
the subgroup $\pi_1(RK) \subset \pi_1(RK/T)$ of the index 2 and
the mapping $RK/T \to K(\I_a,1)$ is well-defined. The image of the
element $c$ is the generator of the cyclic group $\I_a$.) The
composition of the left vertical arrows in the diagram $N_d \to
N_d/T_{N_d} \to K(\I_a,1)$ coincides with the mapping $\mu: N_d
\to K(\I_a,1)$. The pair $(N^{n-2k}, \mu_a)$ is cobordant to a
pair $(N'^{n-2k}, \mu'_a)$, where $N'^{n-2k} = N'^{n-2k}_{cycl}
\cup N'^{n-2k}_d$, and $N'^{n-2k}_{cycl}$ is a closed manifold.
The manifold (with boundary) $N'^{n-2k}_{cycl}$ is the double
covering over the oriented manifold (with boundary)
$N'^{n-2k}_{cycl}/T_{N_d}$. The base of the cover represents a
cycle in  $H_{n-2k}(K(\I_a,1),K(\I_d,1);\Z)$. Therefore the
relative cycle $\mu'_{a,\ast}([N'^{n-2k}_d,\partial N']) \in
H_{n-2k}(K(\I_a,1),K(\I_d,1))$ is trivial and
$$  \left< \mu^{\ast}_a (t) ; [N^{n-2k},\partial N^{n-2k-1}]
\right>=  \left< \mu'^{\ast}_a  (t) ; [N'^{n-2k}_{cycl}] \right>.
$$
The last characteristic number coincides
with
$$ \left< \kappa^{n-k}; [\bar N^{n-k}] \right> = 1. $$
Lemma 8 is proved.

\subsubsection*{A natural stratification of the polyhedron $K$. Proof of Lemma 7}
Let $J$ be the join of lens spaces $(S^7/i)_j$, $j = 1 \dots, r$.
The space $J$ admits a natural stratification defined by the
collection of subjoins $J(k_1, \dots, k_s)$ generated by lenses
with numbers $0< k_1 < \dots < k_s < r$.

The preimage $$p^{-1}(J(k_1, \dots, k_s)) \subset \RP^{n-k},$$ $p:
\RP^{n-k} \to J$, is denoted by $R(k_1, \dots, k_s)$.


A point $$ x \in R(k_1, \dots, k_s) \subset \RP^{n-k}$$ is
determined by the collection of coordinates $$(x_{k_1}, \dots,
x_{k_s}, \lambda)$$ (up to the antipodal transformation of the first $s$ coordinates), where
$x_{k_j} \in S^7_j$, and $\lambda$ is a barycentric coordinate on
the standard $(s-1)$-simplex.

The polyhedron $K$ admits a natural stratification
 $$K(k_1, \dots, k_s), \quad 1 \le s \le r, $$ correspondingly to the stratification of $J$.
 The maximal stratum $K(1, \dots, r)$ is represented by the disjoin union of connected components
 of different types.


Let a point $(x_1, x_2)$ belongs to $K(1, \dots, r)$. Let
 $(x_{1,1}, x_{2,1}, \dots,  x_{1,r}, x_{2,r}, \lambda)$ be the collection of the coordinates of the point.
 The first $2r$ terms of this collection is $r$ ordered pairs of points on the standard sphere $S^7$. The collection
 is defined up to the permutation of the coordinates in the pair and up to the independent antipodal transformation of
 first point or second point in each pair.
 The equivalence class of a collection of the coordinates of a point $(x_1,
 x_2)$ contains $8$ collections.

\subsubsection*{Types of components of the maximal stratum}

Let $x \in K(1, \dots, r)$ be a point with the prescribed pair of
collections of spherical coordinates  $(x_{1,i},x_{2,i})$. The
following possibilities are: the coordinates in the $i$-th pair

 (1) coincide, or (2) are antipodal, or
(3) are related by means of the generator of the $\Z/4$-cyclic
cover.

 This determines a
sequence of $r$ complex numbers $v_i \in \{1,-1,+i,-i\}$, $i=1,
\dots, r$, with respect to (1),(2), or (3). We will call such a
sequence the characteristic. For an arbitrary point in a
prescribed component $K(1, \dots, r)$ the characteristic
 is well-defined up to the multiplication of each term by
$-1$ and this characteristic does not depend on a point on the
component. We shall say that the prescribed component of the
maximal stratum is of the $\I_a$-type ($\I_b$-type) if the
corresponded characteristic contains only $\{+i,-i\}$
($\{+1,-1\}$); the component is of the $\I_d$-type, if the
characteristic contains at least $3$ different values. It is easy
to prove that the restriction of the canonical mapping $\eta :
\Gamma \to K(\D_4,1)$ on a stratum of the $\I_a,\I_b$ or
$\I_d$-type admits (up to homotopy equivalence of the mappings) a
reduction with the target in the subspace $K(\I_a,1)$,
$K(\I_b,1)$, $K(\I_d,1)$ of the space $K(\D_4,1)$ correspondingly.
This reduction (for strata of the $\I_a$ and $\I_b$-types) is
well-defined up to the composition with mapping of the
corresponding classified space given by the conjugation
automorphism $\D_4 \to \D_4$, $x \to (ba)x(ba)^{-1}$, $x \in
\D_4$, $ba \in \I_c$.


\subsubsection*{The resolution space $RK$}

Let us denote by $K_1 \subset K$ the disjoint union of all
singular strata of the length $1$, by  $K_0$ the disjoin union of
maximal strata, and by $K_{reg} \subset K$-- the subpolyhedron
defined as $K_{reg}= K_0 \cup K_1$. The component of the boundary
$K_{reg} \cap Q_{antidiag}$ is denoted by $Q_{reg,antidiag}$ and
the component of the  boundary $K_{reg} \cap Q_{diag}$ is denoted
by $Q_{reg,diag}$. Note that $Q_{reg,antidiag}$
($Q_{reg,antidiag}$) contains only points for which no more then
two numbers in the characteristic are different from $+i$ ($+1$).
Components of the space $K_0$ are divided into 3 classes:
diagonal, antidiagonal and generic class. A component of the
diagonal (antidiagonal) class
 intersects with the diagonal (the
 antidiagonal) by a maximal subcomponent of the boundary.  A regular stratum of the considered type
 contains only points for which the only number  in the characteristic is
different from $+i$ ($+1$).

Let us denote by $\bar K^1$ the 2-sheeted covering space over
$K^1$ with respect to the inclusion $\I_c \subset \D_4$. This
covering coincides with  the canonical double covering over the
polyhedron of self-intersection points.

The space $RK$ is defined from the following diagram:
$$ \bar K^1 \to K^1 \subset K \supset K_{reg}.$$
The space $K_{reg}$ is defined by the gluing by means of the collection of 2-sheeted
coverings over a regular neighborhood of each component of the singular stratum of the
length 1 with respect to the mapping $$U(\bar K_{1})\setminus \bar K_1
\longrightarrow (K_{reg} \setminus K_1).$$


\subsubsection*{The cyclic mapping $\phi: RK \to K(\I_a,1)$}

 The union of all components in the given class is denoted by $K_{0,diag}$,
 $K_{0,antidiag}$,$K_{0,int}$ correspondingly.
The restriction of the mapping $\eta: K \to K(\D_4,1)$ to
$K_{0,diag} \subset K$ ($K_{0,antidiag} \subset K$) is given by
the composition
$$K_{0,diag} \to K(\I_b,1) \subset K(\D_4,1)$$ $$(K_{0,antidiag} \to
K(\I_a,1) \subset K(\D_4,1)).$$

Because of the prescribed boundary condition, the structured
reduction for a subcomponent of $K_{0,diag}$, $K_{0,antidiag}$ is
canonical. The structured reduction for a component of $K_{0,int}$
is non-canonical.

The mapping $$\phi_0: K_{0} \to K(\I_a,1)$$ extends to a mapping
$$\phi: RK \to K(\I_a,1).$$
\[  \]
The composition $$\bar K_{1} \to K_{1}
\stackrel{\eta}{\longrightarrow} K(\D_4,1)$$ admits a natural
reduction with the target
 $$K(\I_c,1)
\subset K(\D_4,1).$$

Let $K_{0,\alpha}$, $K_{0,\beta} \subset K_0$ be two prescribed
components of the same type $\I_b$ (or $\I_a$) with a common
boundary stratum $K_{1,\gamma} \subset K_1$.
 A cyclic mapping
$$\phi_{0,\ast} = \pi_d \circ \eta_{\ast}: K_{0,\ast} \to K(\I_b,1)
\to K(\I_d,1)$$ $$(\eta_{\ast}: K_{0,\ast} \to K(\I_a,1)), \quad
\ast \in \{\alpha, \beta\},$$ where $K_{0,\ast} \to K(\I_b,1)$
($K_{0,\ast} \to K(\I_a,1)$) is well-defined up to a composition
with the mapping $$K(\I_b,1) \to K(\I_b,1)\quad (K(\I_a,1) \to
K(\I_a,1)).$$ The last mapping is induced by the automorphism
$$\D_4 \to \D_4, \quad x \to (ba)x(ba)^{-1},$$ $$ x \in \D_4,
\quad ba \in \I_c.$$

The transfer with respect to the inclusion $\I_c \subset \D_4$
determines a unique mapping $$\eta^!_{\ast}: \bar K_{0,\ast} \to
K(\I_d,1).$$ This proves that the extension $\phi$ with the
prescribed boundary conditions exists.


\subsubsection*{The lift $ res: N_d \to RK$}
Let us consider a generic $PL$-homotopy $$F(\tau): S^{n-k}/i \to
\R^n, \quad \tau \in [0;1]$$ with the boundary conditions
$$F(0)=i \circ \hat p: S^{n-k}/i \to J \subset \R^n.$$  For a given
$\tau \in (0;1]$ the double points of $F(\tau)$ is denoted by
$\hat N(\tau)$. This manifold with boundary is a submanifold of a
quotient of the space $\Gamma/T_{\Gamma} \times \{\tau\}$. The
polyhedron $\cup_{\tau} \hat N(\tau)$, $\tau \in (0,\varepsilon]$
($\varepsilon$ is sufficiently small) is denoted by $\hat
N_{(0;\varepsilon]}$.


Because the mapping $F$ is a $PL$-mapping, the bottom boundary of
$\hat N_{(0;\varepsilon]}$, denoted by $\hat N_0$, is a
$15$-dimensional subpolyhedron in a quotient of the space
$\Gamma/T_{\Gamma} \times \{0\}$. The polyhedron $\hat
N_{(0;\varepsilon]}$ is the base of the 4-sheeted cover
$N_{(0;\varepsilon]} \to \hat N_{(0;\varepsilon]}$, where
$N_{(0;\varepsilon]}$ is the set of  self-intersection points of
the composition (after a small alteration)
$$F(\tau) \circ \pi:
\RP^{n-k} \to S^{n-k}/i \to \R^n.$$

 Because
of the general position arguments the following condition holds:



(1) the polyhedron $N_0$ does not intersect (if $\varepsilon$ is
small enough) with singular strata in $\Gamma$ of the length $\ge
2$ (of codimension $ \ge 16$).
\[  \]
(2) the polyhedron $N_0$ is in  general position with respect to
the stratum of  length 1; in particular, the restriction of
$F(\tau) \vert_{p^{-1}(J^1)}$, $J^1 \subset J$, $\tau \in
(0,\varepsilon]$ to the singular stratum of  length 1  is an
embedding.

A {\it resolution
 mapping}
 $$ res: N_d(\varepsilon) \to RK$$
 with the
prescribed boundary conditions is well-defined from (1),(2). Note
that $diam(U_{antidiag}), diam(U_{diag})$ has to be less then the
distance between $N_0$ and $K_2 \subset K$.

 Proposition 6 is proved.

\section*{Discussion}

\subsubsection*{Conjecture 1}

There exists a 7-dimensional manifold $K^7$ with a normal
$\D_4$-framing $\Xi_K$ in codimension $2^l-8$, $l \ge 4$, such
that the pair $(K^7,\Xi_K)$ has the Steenrod-Hopf invariant 1.


\subsubsection*{Remark}
An arbitrary cyclic $\I_a$-framed manifold $(N^7,\Xi_N)$ in
codimension $2^l-8$ has the trivial Steenrod-Hopf invariant. The
conjectured $\D_4$-framed manifold $(K^7, \Xi_K)$ cannot be
realized as a double-point manifold for a skew-framed immersion
$f: M^{2^{l-1}+3} \looparrowright \R^{2^{l}-1}$.


\subsubsection*{Conjecture 2}

The Main Theorem holds for $n \ge 31$, i.e. for an arbitrary
skew-framed immersion $f: M^{2^{l-1}+7} \looparrowright
\R^{2^{l}-1}$, $l \ge 5$  the Steenrod-Hopf invariant is trivial.

\subsubsection*{Remark}
A proof of the Conjecture 2  could be obtained by means of a
straightforward generalization of the Melikhov mapping. The join
of the standard mappings
$$\RP^7 \to S^7/i \subset \R^{14}$$  is
replaced by the join of several copies of the standard mapping
$$\RP^{3} \to Q^3 \subset \R^4,$$
where $Q^3$ is the quotient of the 3-sphere by the quaternions
group of the order 8 (a homogeneous space) standard space,
$$\RP^3 \to Q^3$$ is the standard 4-sheeted cover, $$Q^3 \subset
\R^4$$ is the Massey embedding. This embedding is explicitly
described in [M], Example 4.

This generalized construction determines the mapping
$$\RP^{4k-1} \to \R^{5k-1}$$ (below the metastable range) with a
cyclic structure. The case $$k=6, \quad 4k-1 \,=\, 23\, =\,
2^4+7,$$
$$5k-1\, =\, 29\, <\, 31\, =\, 2^5-1$$ is required for a generalization of the The Main Theorem 1.
\[  \]

The present paper was started at Postnikov's Seminar in 1996 and
was finished at Prof. A.S.Mishenko Seminar. This paper is
dedicated to the memory of Prof. M.M.Postnikov. The paper was
presented at the M.M.Postnikov Memorial Conference (2007)--
"Algebraic Topology: Old and New". A preliminary version was
presented at the Yu.P.Soloviev Memorial Conference (2005)
"Topology, analysis and applications to mathematical physics".

\end{document}